\documentclass[12pt]{amsart}
\usepackage[latin1]{inputenc}
\usepackage{mathptmx}
\usepackage{amscd}
\theoremstyle{theorem}
\newtheorem{theorem}{Theorem}
\newtheorem{lemma}[theorem]{Lemma}

\theoremstyle{definition}

\let\Bbb=\mathbb
\let\phi=\varphi
\def\NZQ{\Bbb}
\def\CC{{\NZQ C}}
\def\ZZ{{\NZQ Z}}
\def\RR{{\NZQ R}}
\def\pp{{\frak p}}

\def\grade{\operatorname{grade}}

\def\Ass{\operatorname{Ass}}
\def\Min{\operatorname{Min}}

\def\degree{\operatorname{deg}}
\def\lcm{\operatorname{lcm}}
\def\gcd{\operatorname{gcd}}

\let\oldbigwedge\bigwedge
\def\BIGwedge{{\textstyle\oldbigwedge}}
\def\medwedge{{\scriptstyle\oldbigwedge}}
\def\bigwedge{\mathchoice{\BIGwedge}{\BIGwedge}{\medwedge}{}}

\def\pp{{\mathfrak p}}

\let\iso=\cong

\let\epsilon=\varepsilon

\begin{document}
\title{On the coefficients of Hilbert quasipolynomials}

\author{Winfried Bruns and Bogdan Ichim}
\address{Universit\"at Osnabr\"uck, FB Mathematik/Informatik, 49069
Osnabr\"uck, Germany} \email{winfried@math.uos.de}

\address{Universit\"at Osnabr\"uck, FB Mathematik/Informatik, 49069
Osnabr\"uck, Germany, and \phantom{iii} Institute of Mathematics, C.P. 1-764,
70700 Bucharest, Romania} \email{bogdan.ichim@math.uos.de,
bogdan.ichim@imar.ro}

\begin{abstract}
The Hilbert function of a module over a positively graded algebra is
of quasi-polynomial type (Hilbert--Serre). We derive an upper bound
for its grade, i.~e.\ the index from which on its coefficients are
constant. As an application, we give a purely al\-gebraic proof of
an old combinatorial result (due to Ehrhart, McMullen and Stanley).
\end{abstract}

\maketitle

\section{Hilbert quasipolynomials}

Let $K$ be a field, and $R$ a positively graded $K$-algebra, that
is, $R=\bigoplus_{i\ge 0}R_i$ where $R_0=K$ and $R$ is finitely
generated over $K$. We do not assume $R$ to be generated in degree
$1$ -- the generators may be of arbitrarily high degree. The theorem
of Hilbert--Serre describes the Hilbert functions of finitely
generated graded $R$-modules $M$:

\begin{theorem}\label{Serre}
Let $M=\bigoplus_{i\in \ZZ} M_i$ be a finitely generated graded
$R$-module of dimension $d$, $H(M,\_):\ZZ\to\ZZ$ the associated
Hilbert function, and suppose that $r_1,\dots,r_d$ is a homogeneous
system of parameters for $M$.

Then there is a quasi-polynomial $Q_M$ of degree $d-1$, such that
$H(M,n)=Q_M(n)$ for $n\gg 0$. Moreover, the period of $Q_M$ divides
$\lcm(\deg r_1,\dots,r_d)$.
\end{theorem}

The terminology concerning quasipolynomials is explained as follows:
a function $Q:\ZZ\to\CC$ is called a \emph{quasi-polynomial} of
\emph{degree} $u$ if
$$
Q(n)=a_u(n)n^u+a_{u-1}(n)n^{u-1}+\ldots+a_1(n)n+a_{0}(n),
$$
where $a_i:\ZZ\to\CC$ is a periodic function for $i=0,\dots,u$, and
$a_u\neq 0$. The \emph{period} of $Q$ is the smallest  positive
integer $\pi$ such that
$$
a_i(n+m\pi)=a_i(n)
$$
for all $n,m\in \ZZ$ and $i=0,\dots,u$.

For the reader's convenience, we include a short proof the
Hilbert--Serre theorem, or rather its reduction to the classical
theorem of Hilbert. By definition of homogeneous system of
parameters, $M$ is a finitely generated module over
$K[r_1,\dots,r_d]$ (which is isomorphic to a polynomial ring over
$K$). Therefore we may assume that $R=K[r_1,\dots,r_d]$. Let $S$ be
the subalgebra of $R$ generated by its homogeneous elements of
degree $p=\lcm(\deg r_1,\dots,\deg r_d)$. Then it is not hard to see
that $R$ is a finitely generated $S$-module. Therefore $M$ is a
finitely generated $S$-module, too, and $\dim_S M=\dim_R M$. As a
last reduction step, we can replace $R$ by $S$ and assume that $R$
is generated by its elements of degree~$p$.

Then we have the decomposition
$$
M=M^0\oplus\ldots\oplus M^{p-1},\qquad M^{k}=\bigoplus_{i\equiv k\
(p)}M_i,
$$
into finitely generated $R$-modules, and $\dim M=\max_k \dim M^k$.

Let us consider a single module $M^k$. Then we can normalize the
degrees in $R$ dividing them by $p$ and re-grade $M^k$ by giving
degree $(i-k)/p$ to the elements of its degree $i$ component in the
original grading, $i\equiv k\ (p)$. By Hilbert's theorem, the
Hilbert function of $M^k$ re-graded is given by a true polynomial
$P_k(n)$ for $n\gg 0$.

Returning to $M$ we obtain
$$
H(M,n)=P_k((n-k)/p),\qquad n\equiv k\ (p),\ n\gg 0,
$$
and this proves the theorem.

It is clear that any improvement of the theorem depends on the
``coherence'' of the modules $M^k$. The reduction in the proof above
forgets the original module structure to a large extent. Clearly, in
the extreme case in which $R$ is generated by its degree $p$
elements, $M$ is just a direct sum of the independent modules $M^k$.
But if the $M^k$ are sufficiently related, then one can say more on
$Q_M$.

\section{The grade of Hilbert quasipolynomials}

It is a natural question to ask how close $Q_M$ is to being a true
polynomial. The next theorem, which is the main result of this
paper, provides an answer. Following Ehrhart \cite{E}, we let the
\emph{grade} of $Q$ denote the smallest integer $\delta\ge -1$ such
that $a_i(\_)$ constant for all $i>\delta$.

\begin{theorem}\label{main}Let $M=\bigoplus_{i\in \ZZ} M_i$ be a
finitely generated graded $R$-module of dimension $d$, and
$$
Q(n)=a_{d-1}(n)n^{d-1}+a_{d-2}(n)n^{d-2}+\ldots+a_1(n)n+a_{0}(n)
$$
its Hilbert quasi-polynomial with period $\pi$. Let $I$ be the ideal
of $R$ generated by all homogeneous elements $x$ of $R$ such that
$\gcd (\degree x,\pi)=1$.  Then
$$
\grade Q<\dim M/IM.
$$
\end{theorem}

The theorem will be proved by an induction based on the following
lemma, in which, as usual, $(0:x)_M=\{u\in M: xu=0\}$.

\begin{lemma}
With the notation of the theorem, if $\dim M/IM<\dim M$, then there
is a homogeneous $x\in I$ with $\gcd (\degree x,\pi)=1$, such that
\begin{enumerate}
\item[\rm{(a)}] $\dim M/xM=\dim M-1$,
\item[\rm{(b)}] $\dim (0:x)_M\le \dim M-1$.
\end{enumerate}
\end{lemma}

\begin{proof}
Let $D(M)=\{\pp\in V(M),\dim A/\pp=\dim M\}=\{\pp_1,\ldots,\pp_r\}$.
Clearly $I\not\subset \pp_i$ for $i=1,\dots,r$. By prime avoidance,
we conclude that $I\not\subset\pp_1\cup\dots\cup\pp_r$. By induction
on $r$, we show that
$$
S=\{x\in I,\ x \text{ homogeneous},\ \gcd (\degree
x,\pi)=1\}\not\subset\bigcup_{i=1}^r \pp_i.
$$
This is clear for $r=1$. For $1\le j\le r$, we may assume by
induction that
$$
S\not\subset\bigcup_{i=1,i\ne j}^r \pp_i.
$$
Assume that $S\subset\pp_1\cup\dots\cup\pp_r$. Then for each
$j=1,\dots,r$  there is $x_j\in S$ such that
$$
x_j\in \pp_j\setminus \big(\bigcup_{i=1,i\ne j}^r \pp_i\big).
$$
Let $\degree x_1=\alpha$ and $\degree x_2\cdots x_r=\beta$. Then
$x=x_1^{\lcm(\alpha,\beta)/\alpha}+(x_2\ldots
x_r)^{\lcm(\alpha,\beta)/\beta}\in S$, since it is homogeneous, and
$\gcd(\lcm(\alpha,\beta),\pi)=1$. Now
$$
x_1\in \pp_1\setminus \big(\bigcup_{i=2}^r
\pp_i\big)\quad\text{and}\quad x_2\ldots x_r\in\big(\bigcap_{i=2}^r
\pp_i\big)\setminus \pp_1\quad\text{implies}\quad
x\not\in\bigcup_{i=1}^r \pp_i,
$$
a contradiction.

Let $x\in S\setminus (\pp_1\cup\dots\cup\pp_r)$. Then $\dim
M/xM=\dim M-1$. Moreover every prime ideal in then support of
$((0:x)_M$ is in the support of $M/xM$. Thus $\dim (0:x)_M\le \dim
M-1$.
\end{proof}

\begin{proof}[Proof of Theorem \ref{main}]
We prove by induction on  $\dim M=d$ that $\dim M/IM\le\gamma$
implies $a_j(\_)$ constant for all $j\ge \gamma$. This is clear if
$d\le\gamma$ (then $j\ge \gamma$ implies $a_j(\_)=0$), so we may
assume $d>\gamma$. Let $x$ be as in the lemma, and $g=\degree x$.

Set $M'=M/xM$ and $M''=(0:x)_{M}$. Then $M'/IM'\iso M/IM$ and
certainly $\dim M''/IM''\le \gamma$. Since $\dim M',\dim M'' <\dim
M$, we may assume by induction that
$H(M/xM,n)$ and $H((0:x)_{M},n)$, $n\gg0$, are quasipolynomials of
grade $< \gamma$.

The exact sequence
$$
\CD 0@>>>(0:x)_{M}(-g)@>>>M(-g)@>x>>M@>>>M/xM@>>>0
\endCD
$$
gives the equation
$$
H(M,n)-H(M,n-g)=H(M/xM,n)-H((0:x)_{M},n-g).
$$
For a quasipolynomial $Q$ it is easy to see that $Q(n-g)$ has the
same grade as $Q$. Therefore the right hand side in the previous
equation is a quasipolynomial of grade $<\gamma$ for $n\gg0$, and so
this holds for the left hand side, too. So it remains only to apply
the following lemma.
\end{proof}

\begin{lemma}
Let $Q(n)=\sum a_k(n)n^k$ be a quasipolynomial. If $Q(n)-Q(n-g)$ is
of grade $<\gamma$ for some $g$ coprime to the period $\pi$ of $Q$,
then $\grade Q<\gamma$.
\end{lemma}

\begin{proof}
Let $u=\deg Q$ and let us first compare the leading coefficients. We
can assume $\gamma\le u$. Then one has $a_u(n)-a_u(n-g)=C$ for some
constant $C$ and all $n$, and so $a_u(n)-a_u(n-\pi g)=\pi C$. Since
$\pi$ is the period, we conclude that $C=0$, and $a_u(n)=a_u(n-g)$.
But $g$ is coprime to $\pi$, and it follows that $a_u$ is constant.

The descending induction being started, one argues as follows for
the lower coefficients. Suppose that $k\ge\gamma$. Then
$a_k(n)-a_k(n-g)$ is a polynomial in the coefficients $a_j$ for
$j>k$ and $g$. Since the higher coefficients are constant by
induction, it follows that $a_k(n)-a_k(n-g)$ is constant, too, and
the rest of the argument is as above.
\end{proof}

\section{An Application to Rational Polytopes}

In this section we shall give a purely algebraic proof of an old
theorem, which was conjectured by Ehrhart (\cite{E} p.\ 53), and
proved independently by  McMullen (see \cite{M}) and Stanley
(\cite{S}, Theorem 2.8):

\begin{theorem}\label{application} Let $P$ be a $d$-dimensional rational
convex polytope in $\RR^m$, and let the Ehrhart quasi-polynomial of
$P$ be
$$
E_P(n)=a_{d}(n)n^{d}+a_{d-1}(n)n^{d-1}+\ldots+a_1(n)n+a_{0}(n).
$$
Suppose that for some $\delta$ the affine span of every
$\delta$-dimensional face of $P$ contains a point with integer
coordinates. Then $\grade E_P< \delta$.
\end{theorem}

\begin{proof}
We choose a field $K$ and let $R$ be the Ehrhart ring of $P$. It is
the vector subspace of $K[X_1^{\pm1},\dots,X_m^{\pm1},T]$ spanned by
all Laurent monomials $X^aT^n=X_1^{a_1}\cdots X_m^{a_m}T^n$ where
$a=(a_1,\dots,a_m)\in nP$, $n\in\ZZ$, $n\ge 0$. By Gordan's lemma it
follows easily that $R$ is a finitely generated, positively graded
$K$-algebra, where we use the exponent of $T$ as the degree of a
monomial. The Ehrhart function of $P$ is just the Hilbert function
of $R$. (See Chapter 6 of \cite{BH} for more information.)

Let $\pi$ be the period of $E_P$, and $F$  a $\delta$-dimensional
face  of $P$. Since the affine span of $F$ contains a point with
integer coordinates, $nF$ contains a point with integer coordinates
for all $n\gg 0$. We chose $n$ big enough so that $nF$ contains a
point $m_F$ with integer coordinates for every $\delta$-dimensional
face $F$, and $\gcd(n,\pi)=1$.

Now let $J\subset R$ be the ideal generated by the monomials
$X^{m_F}T^n$. If $\dim R/J\le\delta$, then we are done by Theorem
\ref{main} because the ideal $I$ in Theorem \ref{main} contains $J$.

Since $J$ is a monomial ideal, $\Ass_R R/J$ consists of monomial
prime ideals. In particular, $\Min_R R/J$ consists of monomial prime
ideals. By theorem 6.1.7 of \cite{BH}, for each $\pp\in\Min_R R/J$,
there is a face $G_\pp$ of $P$, such that $\pp$ is generated by all
monomials outside the cone associated with $G_\pp$. One has $\dim
R/\pp=\dim G_\pp+1$. Since $J\subset \pp$, it follows that $\dim
G_\pp\le \delta-1$. So $\dim R/J=\max\{\dim R/\pp, \pp\in \Min_R
R/J\}\le(\delta-1)+1=\delta$.
\end{proof}

\end{document}